\documentclass{csri21}
\newcommand{\state}{\boldsymbol u}

\newcommand{\resid}{\mathbf{r}}
\newcommand{\spacetimeResid}{\overrightarrow{\mathbf{r}}}

\newcommand{\params}{\boldsymbol \mu}
\newcommand{\paramDomain}{\mathcal{D}}

\newcommand{\RR}[1]{\mathbb{R}^{#1}}

\newcommand{\basis}{\boldsymbol \Phi}
\newcommand{\velocity}{\boldsymbol f}

\newcommand{\romDim}{K}

\newcommand{\density}{\rho}

\newcommand{\expectation}{\mathbb{E}}
\newcommand{\variance}{\text{var}}
\newcommand{\bz}{\boldsymbol 0}

\newcommand{\timeStep}{\Delta t}
\newcommand{\mA}{\mathbf{A}}
\newcommand{\vb}{\mathbf{b}}

\newcommand{\vState}{\mathbf{u}}
\newcommand{\spacetimeState}{\overrightarrow{\mathbf{u}}}
\newcommand{\reducedState}{\hat{\mathbf{u}}}
\newcommand{\mState}{\mathbf{U}}
\newcommand{\identity}{\mathbf{I}}

\newcommand{\vm}{\mathbf{m}}

\usepackage{xcolor}
\usepackage{tikz}
\usetikzlibrary{shapes,arrows}
\tikzstyle{block} = [rectangle, draw, fill=blue!20, text width=10em, text centered, rounded corners, minimum height=4em]
\tikzstyle{line} = [draw, -latex']
\usepackage{amsmath}
\usepackage{cleveref}
\usepackage{booktabs}

\usepackage{amsfonts,amsmath,graphicx,subfigure}

\addcontentsline{toc}{chapter}{The Quantum Mechanics of Chocolate Pudding\\
{\em S.S.\ Student and M.M.\ Mentor}}

\title{Space-time reduced-order modeling for uncertainty quantification}

\author{Ruhui Jin \thanks{Department of Mathematics, University of Texas at Austin, rhjin@math.utexas.edu} \and Francesco Rizzi \thanks{NexGen Analytics, fr.francescorizzi@gmail.com} \and Eric Parish \thanks{Sandia National Laboratories, ejparis@sandia.gov}}

\begin{document}

\maketitle

\begin{abstract}
This work focuses on the space-time reduced-order modeling (ROM) method for solving large-scale uncertainty quantification (UQ) problems with multiple random coefficients.
In contrast with the traditional space ROM approach, which performs dimension reduction in the spatial dimension, the space-time ROM approach performs dimension reduction on both the spatial and temporal domains, and thus enables accurate approximate solutions at a low cost. 
We incorporate the space-time ROM strategy with various classical stochastic UQ propagation methods such as stochastic Galerkin and Monte Carlo.
Numerical results demonstrate that our methodology has significant computational advantages compared to state-of-the-art ROM approaches. By testing the approximation errors, we show that there is no obvious loss of simulation accuracy for space-time ROM given its high computational efficiency.
\end{abstract}

\section{Introduction}
Quantifying uncertainties in physical systems plays an important role in numerous fields, including climate modeling \cite{GM12}, hypersonic aerodynamics \cite{BDW04, LFGCCM14} and quantum mechanics \cite{OCHHR17}. It has long been a computational challenge to model and simulate large-scale dynamical and control systems with high-dimensional parametric uncertainties.
Researchers have been developing model reduction methods \cite{LBS04, BOCW17} to tackle this computational bottleneck. By building and working with a reduced-order model (ROM) as a qualified approximation to the full-oder model (FOM), the overall computational complexity is reduced significantly. 

Current ROM studies mostly consider the space ROM method and focus on only spatial dimension reduction but maintain the full dimensionality of the temporal domain. As a result, space ROMs can have limited computational savings for unsteady problems characterized by, e.g., small required time steps or long simulation horizons.
Regarding the UQ approach, on one hand, Monte Carlo (MC) is by far the most popular method applied in the ROM workflow \cite{PGP13, UHLM15} due to its reliability and implementation simplicity. On the other hand, other types of UQ propagation methods, for example, the stochastic Galerkin (SG) technique \cite{BTG04, LCE18} based on polynomial chaos expansion, have advantages of good spectral accuracy and convergence over the classical MC.

In this work, we study the space-time ROM method \cite{BBH18, CC19} constructed via Galerkin projection and space-time proper orthogonal decomposition. This novel approach is considered a variation of the space ROM. Its implementation simply stacks the space and time dimensions to achieve the model reduction by finding a lower-dimensional representation for both spatial and temporal domains. The method simultaneously approximates a large-scale PDE model for all points in space and time within a much faster computing time compared to the commonly used space ROM method. Additionally, the space-time ROM approach is often equipped with more favorable error bounds and stability properties than space ROMs~\cite{CC19}.

We apply various UQ propagation techniques such as Monte Carlo and stochastic Galerkin in the space-time ROM framework. We test the described methodology on advection-diffusion PDE problems with multi-dimensional parametric uncertainties.
Our numerical results show that the space-time approach can result in huge computational speed-ups while maintaining accurate approximated solutions.

The main contributions of this work are:
\begin{enumerate}
\item
We study the space-time ROM method and combine it with the well established stochastic Galerkin (SG) technique. By constructing polynomial basis on the reduced space-time domain, we demonstrate that the space-time ROM with SG has good accuracy and a faster computing time as compared to traditional ROM techniques (e.g., the space ROM with MC sampling).
\item
We implement the proposed computational scheme for one and two-dimen-sional advection-diffusion-reaction PDE problems.
\item
We provide thorough numerical assessments for the space-time ROM with respect to the computational time, approximation errors and convergence property given increasing number of samples (MC) and polynomial degrees (SG). 
\end{enumerate}

\section{Mathematical background}
\subsection{Full-order model}
We consider the numerical solution to the parametrized dynamical system:
\begin{equation}
\label{problem form}
\dot{\state}(t, \params) = \velocity(\state(t, \params), t, \params)),\quad \state(0,\params) = \state_{\bz}(\params)
\end{equation}
where
\begin{enumerate}
\item
$\params\in \paramDomain \subset \RR{N_{\params}}$ denotes uncertain parameters;
\item
$\state: [0,T] \times \paramDomain \to \RR{N_s}$ is the time-dependent, parametrized state as the solution to problem \eqref{problem form};
\item
$\velocity: \RR{N_s}\times [0,T]\times  \paramDomain \to \RR{N_s}$ is the velocity;
\item
$\state_{\bz} \in \RR{N_s}$ is the initial state.
\end{enumerate}

We aim to understand how the system state $\state(t, \params)$ responds as a function of time $t$ and uncertain parameters $\params.$ To this end, we apply numerical simulation techniques to solve the UQ problem \eqref{problem form}.

We now introduce the time-discretized form of the main problem \eqref{problem form}. In particular, we discretize the temporal domain $[0, T]$ into $N_t$ time instances characterized by $t^n = n\, \timeStep$ where $\timeStep$ denotes the time step. For example, one classical time-discretization method is the Crank-Nicolson method which yields a sequence of discrete solutions $\vState^n(\params) \approx \state(t^n, \params) \in\RR{N_s}$ as the implicit solution to the system of equations at each time step $n = 1, \dots, N_t$:
\begin{equation}
\label{resid at time}
\begin{array}{ll}
\resid^n(\vState^n, \vState^{n-1}, \params)&\displaystyle: \RR{N_s}\otimes\RR{N_s}\otimes\RR{N_{\params}}\to \RR{N_s}  \\
&\displaystyle:= \frac{\vState^n - \vState^{n-1}}{\timeStep} - \frac{1}{2}\, \left(\velocity(\vState^n, t^n, \params) - \velocity(\vState^{n-1}, t^{n-1}, \params)\right)
\end{array}
\end{equation}
with initial condition $\vState^0 = \state_{\bz}(\params).$ Note the parametric dependence of the state has been suppressed in the above for simplicity. Thus a discrete representation of the FOM system is
\begin{equation*}
[\vState^1(\params), \vState^2(\params), \dots, \vState^{N_t}(\params)] \in \RR{N_s} \otimes \RR{N_t}.
\end{equation*}

\subsection{Projection-based model reduction}
The FOM solving process is computationally expensive in practice when the spatial dimension $N_s$ and temporal dimension $N_t$ are large. The reduced-order modeling technique is proposed to overcome this computational challenge. It follows an offline-online paradigm. Please see the workflow \Cref{workflow} below.

In the offline phase, we sample and plug in a certain number of uncertain parameter instances into the full-order model and solve the system accordingly.
The obtained sample solutions are collected to form a snapshot matrix. We then identify a low-dimensional subspace by performing proper orthogonal decomposition (POD) for the state snapshots.
The governing equation \eqref{problem form} is projected onto this trial subspace to create a reduced-order model.
The result of this process is a reduced-order model which can be solved more efficiently.

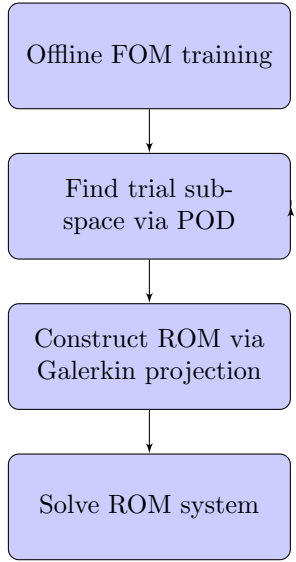
\begin{figure}[tbhp]
\centering
\begin{tikzpicture}[node distance = 2cm, auto]
    \node [block] (offline) {Offline FOM training};
    \node [block, below of=offline] (trial subspace) {Find trial subspace via POD};
    \node [block, below of=trial subspace] (projection) {Construct ROM via Galerkin projection};
    \node [block, below of=projection] (ROM) {Solve ROM system};
    \path [line] (offline) -- (trial subspace);
    \path [line] (trial subspace) -- (trial subspace);
    \path [line] (trial subspace) -- (projection);
    \path [line] (projection) -- (ROM);
\end{tikzpicture}
\caption{ROM workflow}
\label{workflow}
\end{figure}

\subsection{Trial subspace and POD}
\subsubsection{Spatial trial subspace}
Suppose in the offline training procedure, we obtain a collection of snapshot solutions for $N_{\text{train}}$ randomly drawn parameter instances :
\[
\mState_{\text{train}} = \left[\vState^1(\params_1), \dots,\vState^{N_t}(\params_1), \dots, \vState^1(\params_{N_{\text{train}}}),\dots,\vState^{N_t}(\params_{N_{\text{train}}}) \right] \in \RR{N_s \times (N_t\, N_{\text{train}})}.
\]

The proper orthogonal decomposition method identifies a lower-dimensional trial subspace represented by an orthonormal matrix $\basis$  from the above training solution set. In particular, we consider the optimization problem in a least squares sense:
\begin{equation}
\label{POD}
\arg\min_{\basis \in \RR{N_s \times \romDim}} \|\basis\basis^\top \mState_{\text{train}}-\mState_{\text{train}}\|_2^2, \quad \text{subject~to}~\basis^\top\basis = \identity_{\romDim},
\end{equation}
where $\romDim$ is the subspace dimension. For the choice of $\romDim,$ we set a relative energy tolerance threshold $e_{\text{tol}}$ and compute $\romDim$ such that the selected basis $\basis \in \RR{N_s \times K}$ preserves the amount of energy for the training solution set $\mState_{\text{train}}$ that exceeds the threshold. 

To solve \eqref{POD}, we compute the singular value decomposition of $\mState_{\text{train}}$:
$
\mathbf{L}, \mathbf{s},\underline{\hspace{0.25cm}} = \text{SVD}(\mState_{\text{train}}).
$
The subspace dimension $\romDim$ is determined by 
\[
\romDim := \arg\min_{\romDim \in \mathbb{N}} \left\vert e_{\text{tol}}-\frac{\sum_{i=1}^{\romDim} s_i^2}{\|\mathbf{s}\|_2^2}\right\vert.
\]
We then select the first $\romDim$ columns of the left singular vectors $\mathbf{L}$ to form the basis $\basis \in \RR{N_s \times K}$.

\subsubsection{Space-time trial subspace}
In the space-time formulation, instead of expressing residuals of all time steps as in \eqref{resid at time}, we formulate the residual in just one system:
\begin{equation}
\label{spacetime resid}
\spacetimeResid(\spacetimeState(\params), t, \params) = {\bf 0} \in \RR{N_sN_t},
\end{equation}
where we concatenate solutions and residuals of all time steps along one dimension, i.e. $\spacetimeState(\params) = [\vState^{1}(\params)^\top, \dots,\vState^{N_t}(\params)^\top]^\top \in \RR{N_sN_t}.$ 

Identifying a space-time trial subspace is rather similar to identifying the spatial subspace. To be more specific, given a collection of space-time training solutions
\begin{equation*}
\overrightarrow{\mState}_{\text{train}} = [\spacetimeState(\params_1), \dots, \spacetimeState(\params_{\text{train}})] 
\in \RR{(N_sN_t) \times N_{\text{train}}},
\end{equation*}
we apply POD to find a lower-dimensional subspace that captures most of the energy of the above solution set.

\subsection{Galerkin projection}
After identifying a basis $\basis$ of the training solution set, we apply Galerkin projection to construct a reduced-order model. We denote the approximated low-dimensional solution $\reducedState(\params) \in \RR{K}.$ By the assumption of $\vState(\params) \approx \basis \reducedState(\params),$ we impose the residual of the full-order model to be orthogonal to the basis:
\begin{equation}
\label{reduced model}
\basis^\top \resid(\basis\reducedState(\params), t, \params) = {\bf 0} \in \RR{K}
\end{equation}
and solve the above reduced system for $\reducedState(\params)$. 

Note that the above general forms of state solution $\reducedState$ and residual $\resid$ in \eqref{reduced model} can be replaced by $\overrightarrow{\hat{\vState}}, \overrightarrow{{\resid}}$ for the space-time ROM approach.

\section{Uncertainty quantification methods}
In the above section, we introduced the spatial-Galerkin and space-time-Galerkin ROMs to reduce the computational cost associated with solving the forward model. This section details how these ROMs can be combined with several classical UQ propagation methods to solve the underlying UQ problem. In particular, we consider Monte Carlo (MC) sampling and stochastic Galerkin projection. 
\subsection{Monte Carlo sampling}
The MC methodology simply follows as: 
\begin{enumerate}
\item
draw samples of random parameters from certain probability distributions;
\item
solve the system \eqref{problem form} based on these parameter instances;
\item
compute quantities of interest (e.g., mean, variance) from the ensemble of solutions.
\end{enumerate}

\subsection{Stochastic Galerkin}
\label{SG}
We first provide some background of polynomial chaos expansion (PCE) which the stochastic Galerkin approach is built upon. We consider a parametrized linear system
\begin{equation}
\label{parametrized system}
\mA(\params)\vState(\params) = \vb(\params) \in \RR{N},
\end{equation}
where the linear operators $\mA:\paramDomain \to \RR{N}$ and $\vb:  \paramDomain\to \RR{N}$ are constructed correspondingly from the residuals of spatial domain \eqref{resid at time} $(N = N_s)$ or space-time domain \eqref{spacetime resid} $(N = N_sN_t)$, with initial condition given by $\vState^0 = \state_{\bz}(\params) \in \RR{N_s}$. 

The idea is to approximate the numerical solution function $\vState(\cdot): \paramDomain \to \RR{N}$ by using a spectral approximation that lies in the span of a finite set of polynomials $\{\psi_j(\cdot)\}_{j=1}^{N_\psi} \subset L^2(\paramDomain)$. The mathematical formulation is as follows:
\begin{equation}
\label{PCE approx}
\vState(\params)\approx \tilde{\vState}(\params)=\left(\pmb{\psi}(\params)^\top \otimes \identity_N\right)\,\vm,
\end{equation}
where $\pmb{\psi}(\params) = [\psi_1(\params), \dots, \psi_{N_{\psi}}(\params)]^\top \in \RR{N_{\psi}}$ denotes the collection of the polynomial basis and $\vm \in \RR{NN_{\psi}}$ is the coefficient vector .

Thus we can define the approximation residual with respect to polynomial coefficients and uncertainty parameters,
\begin{equation}
\label{PCE resid}
\begin{array}{ll}
\resid(\vm,\params)&: \RR{NN_{\psi}}\otimes \RR{N_{\params}} \to\RR{N} \\
&:= \mA(\params)\,\left(\pmb{\psi}(\params)^\top \otimes \identity_N\right)\, \vm-\vb(\params)\\
&= \left(\pmb{\psi}(\params)^\top \otimes \mA(\params)\right)\, \vm- \vb(\params).
\end{array}
\end{equation}
We compute the unknown coefficients $\vm$ by the residual formula.

We now formally introduce the stochastic Galerkin approach to solve for the coefficients $\vm$ in \eqref{PCE resid}. Given a density function $\density$ for the probability space $\paramDomain,$ we define the inner product:
\begin{equation}
\label{inner product}
\langle g(\params), h(\params)\rangle_{\density} = \int_{\paramDomain} g(\params)\, h(\params)\, \density(\params) d\params,
\end{equation}
where $g, h \in L^2(\paramDomain)$ are functions. The expectation of a function $g$ is given by:
\begin{equation}
\label{sg expectation}
\expectation[g] =  \int_{\paramDomain} g(\params)\, \density(\params) d\params.
\end{equation}

To solve for the coefficients in \eqref{PCE resid}, the stochastic Galerkin method asks to impose orthogonality on the residual of the system \eqref{PCE resid} with respect to the inner product $\langle \cdot, \cdot\rangle$ \eqref{inner product}. That is to say, we restrict the residual to be orthogonal to the polynomial bases, i.e.
\begin{equation*}
\langle\psi_j, r_i(\vm)\rangle_{\density} = \expectation[\psi_j\, r_i(\vm)] = 0,
\end{equation*}
for all residual dimensions $i = 1, \dots, N,$ and stochastic dimensions $j = 1, \dots, N_{\psi}.$ An alternative vector expression is
\begin{equation*}
\expectation[\pmb{\psi}\otimes \resid(\vm)] = {\bf 0}\in  \RR{N\,N_{\psi}}.
\end{equation*}
From the PCE residual formula given in \eqref{PCE resid}, we end up with solving 
\begin{equation}
\expectation[\pmb{\psi\psi}^\top \otimes \mA]\,\vm = \expectation[\pmb{\psi}\otimes \vb]
\end{equation}
and form an approximating function solution $\tilde{\vState}(\cdot)$ by \eqref{PCE approx} at each time step $t = 1, \dots, N_t$. 

Similarly, for the stochastic Galerkin ROM solution $\reducedState(\cdot): \paramDomain \to \RR{K}$, based on the Galerkin projection approach shown in \eqref{reduced model}, we solve coefficients $\hat{\vm} \in \RR{K\,N_{\psi}}$ in a reduced system:
\begin{equation}
\label{sg rom eqn}
\expectation[\pmb{\psi\psi}^\top \otimes \basis^\top\mA\basis]\,\hat{\vm} = \expectation[\pmb{\psi}\otimes \basis^\top\vb]
\end{equation}
and formulate $\reducedState(\cdot)$ by $\left(\pmb{\psi}(\params)^\top \otimes \identity_K\right)\,\hat{\vm}$.

We would like to remark that, for SG space ROM method, we solve \eqref{sg rom eqn} for function $\reducedState(\cdot)$ at each time step, while in the SG space-time ROM scheme, we solve \eqref{sg rom eqn} and obtain function solution $\overrightarrow{\hat{\vState}}(\cdot)$ for all time steps at once. We additionally remark that both the SG space ROM and SG space-time ROM methods result in significantly smaller solution vectors than the standard SG approach and are thus significantly more computationally tractable. 

\section{Numerical experiment I: 1D parametrized advection-diffusion problem}
We focus on the numerical solution of the one-dimensional parametrized advection -diffusion problem with initial and boundary conditions:
\begin{equation}
\label{1d adv-diff}
\left\{
\begin{array}{l}
\displaystyle \frac{\partial \state}{\partial t} + c\, \frac{\partial \state}{\partial x} = \nu\, \frac{\partial^2 u}{\partial x^2},\\\\
\state(0, t, (c, \nu)) = 0,~~\forall~t \in [0,1], (c, \nu) \in \paramDomain\\\\
\state(x, 0, (c, \nu)) = 0,~~\forall~x \in [0,1], (c, \nu) \in \paramDomain
\end{array}
\right.
\end{equation}
where the state $\state: [0,1] \times [0,1] \times \paramDomain \to \RR{}$. Here, the wave speed $c \sim \mathcal{N}(1, 0.15)$ and the diffusion coefficient $\nu \sim U[0.01, 0.02]$. We set the initial condition to be $\reducedState^{0}=\mathbf{0}$ for the space ROM, which exactly enforces the homogeneous initial conditions of the original problem \eqref{1d adv-diff}. The homogeneous initial conditions are trivially satisfied for the space-time ROM. It is noted that more complex initial conditions can be handled by building an affine trial subspace centered about the initial conditions. 

We apply the backward difference scheme for spatial discretization with $N_s = 255$ spatial degrees of freedom. For time discretization, we employ the Crank-Nicolson method with uniform time step $\Delta t = 0.001$ and implicitly solve results for all $N_t = 1/\Delta t = 1000$ time instances. For consistency, we use the same time step in the online phase of both FOM and space ROM approaches and solve the system iteratively at each time step. For the space-time ROM method, we directly obtain the solution at all time steps.

In terms of the ROM solving workflow --- see \Cref{workflow}, we set the number of training samples to be 20 in the offline phase and select a trial subspace that captures at least $99.9999\%$ of the energy in the original snapshot solution set. 

\subsection{Numerical results}
We first consider ROMs equipped with the Monte Carlo sampling approach for the UQ problem \eqref{1d adv-diff}. In what follows, we discuss the computational efficiencies of space ROM and space-time ROM methods in comparison with the FOM solutions. The computation time of a ROM method is calculated as the total running time of:
\begin{equation}
\label{time count}
\text{finding~trial~subspace}+\text{building~ROM~system}+\text{solving~ROM~system}.
\end{equation}

We test FOM and ROM methods on $10,000$ MC samples and give a detailed time assessment for space and space-time ROMs in \Cref{1d_mc_rom_time}. From the table, we can see that the actual solving time (the third column) of the space-time ROM method is $8000$ times lower than the time of space ROM, which shows the computational advantages of the space-time ROM approach. 
\begin{table}[!htb]
{\footnotesize
  \caption{{\bf Computation time for each step in the workflow (unit: second).}}
\label{1d_mc_rom_time}
\begin{center}
  \begin{tabular}{cccc} 
  \toprule[1pt]
 & \bf find trial subspace & \bf build ROM & \bf solve ROM \\[1mm]
\bf space ROM & $0.606$ & $0.005$ & $\bf \color{red}  245.507$\\[1.5mm]
\bf space-time ROM  & 0.197 &$0.370$ & $\bf \color{blue}  0.323$\\ 
 \bottomrule[1pt]
  \end{tabular}
\end{center}
}
\end{table}

\begin{figure}[!htb]
\centering
\includegraphics[width=0.65\textwidth]{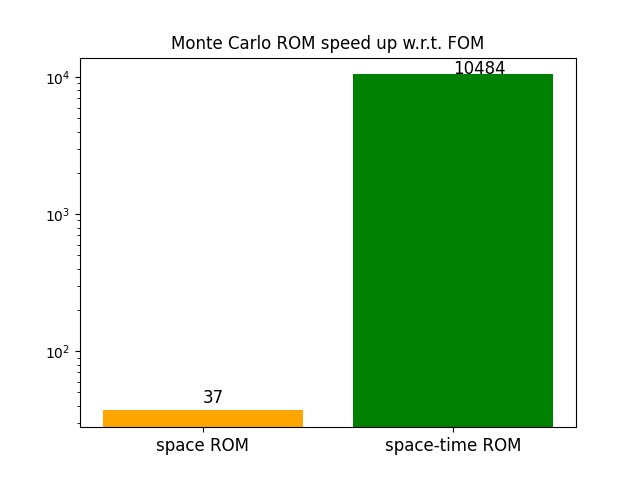}
\caption{Monte Carlo ROM speed-up.}
\label{1d_mc_speedup}
\end{figure}

To evaluate ROM methods' efficiency, we employ the speed-up metric, i.e., the result of FOM's computation time divided by ROM's computation time given the same number of Monte Carlo samples. These time results are obtained by the same experiment for Table 1. It is shown in \Cref{1d_mc_speedup} that the MC space-time ROM method has much greater (3000 times) speed-up compared to the space ROM method.

We proceed to investigate the convergence of the space and space-time ROM methods. Given a certain number of Monte Carlo samples, we plot the relative errors of solution mean and variance respectively for FOM, space ROM and space-time ROM methods. The relative errors for solution mean and variance are defined as follows:
\begin{equation*}
\frac{\left\|\expectation[\vState] - \text{sample~mean}[\basis\reducedState]\right\|_2}{\left\|\expectation[\vState]\right\|_2}, \quad \frac{\left\|\variance[\vState] - \text{sample~variance}[\basis\reducedState]\right\|_2}{\left\|\variance[\vState]\right\|_2},
\end{equation*}
where $\expectation[\vState], \variance[\vState]$ are the expectation and variance of true solution obtained by 12 million FOM MC samples. $\text{sample~mean}[\basis\reducedState]$ and $\text{sample~variance}[\basis\reducedState]$ denote the empirical evaluation for approximated solutions by MC sampling. Note that in the FOM setting, the projection basis $\basis$ is simply the identity matrix. 

\begin{figure}[!htb]
\centering
\subfigure[solution mean]{\includegraphics[width=0.48\textwidth]{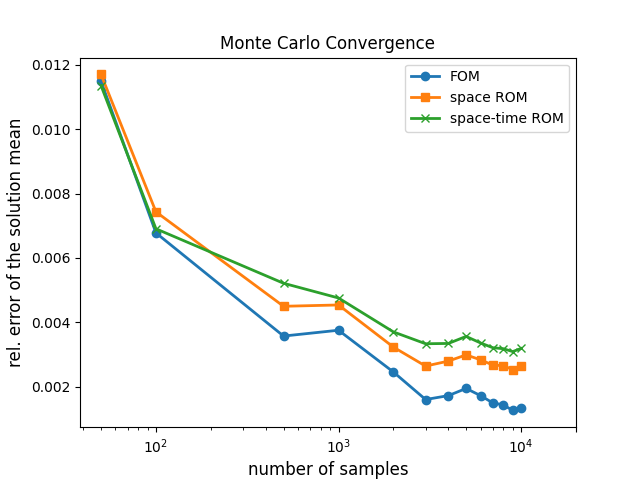}}~~
\subfigure[solution variance]{\includegraphics[width=0.48\textwidth]{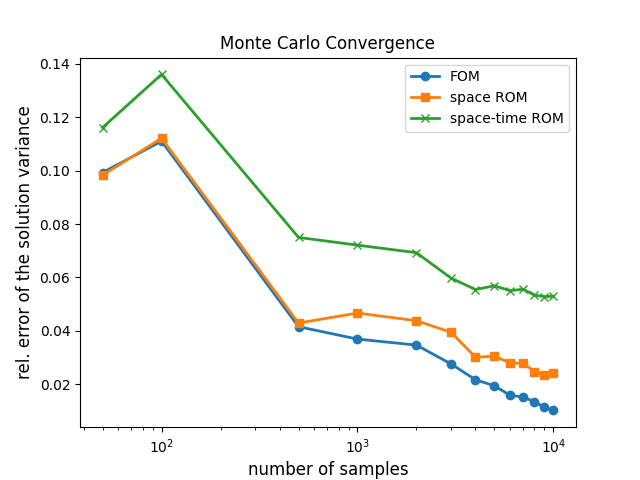}}
\caption{Monte Carlo convergence: relative error versus number of samples.}
\label{1d_mc_convergence}
\end{figure}

Regarding the reproducibility of the experiment, we set a random seed and draw $10,000$ pairs of random samples to run the test. The errors shown in \Cref{1d_mc_convergence} are averaged results from 5 different repetitions of the experiment.

In \Cref{1d_mc_convergence}, all three methods converge well as the number of MC samples increases. The MC FOM method achieves the lowest solution errors, while the MC space-time ROM method in general has the largest errors, especially for the solution variance. We emphasize that the accuracy of both the space and space-time ROMs can be improved by using more basis vectors. Moreover, from the observation on the curves' tendency, the stability of the three methods are rather similar. These numerical results imply space-time approach in general owns good accuracy and stability properties. 

Another important UQ propagation method is stochastic Galerkin based on polynomial chaos expansion (PCE) --- see \Cref{SG} for details. Following the same layout in the Monte Carlo case, we are interested in the computational speed-ups and convergency properties of stochastic Galerkin ROM approaches. 

\begin{figure}[!htb]
\centering
\includegraphics[width=0.65\textwidth]{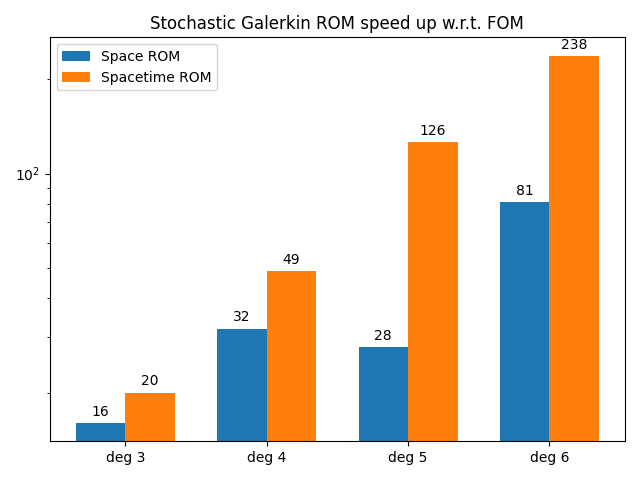}
\caption{stochastic Galerkin ROM speed-up. }
\label{1d_sg_speedup}
\end{figure}

In order to demonstrate the computational efficiency, we plot the space and space-time ROM methods' speed-ups depending on various SG approximation polynomial degrees in \Cref{1d_sg_speedup}. The speed-up gradually increases as the polynomial degree grows. Moreover, the space-time ROM method achieves greater speed-ups than the space ROM.\footnote{Computational times are reported from one run only. We expect that averaging over multiple runs will smooth the observed trends; in particular, the space ROM at polynomial degree 5.}

\begin{figure}[!htb]
\centering
\subfigure[solution mean]{\includegraphics[width=0.48\textwidth]{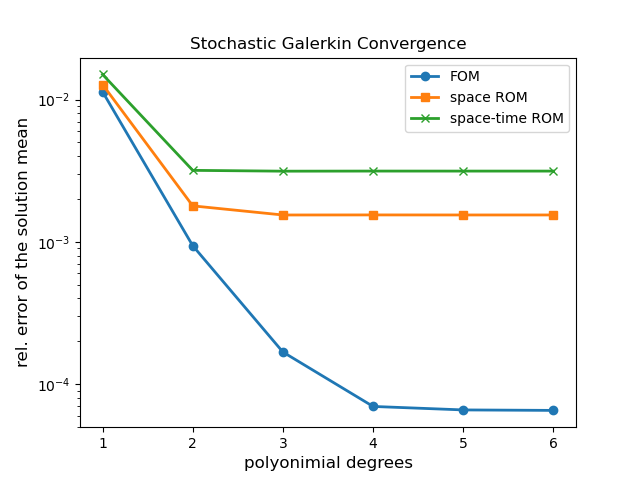}}~~
\subfigure[solution variance]{\includegraphics[width=0.48\textwidth]{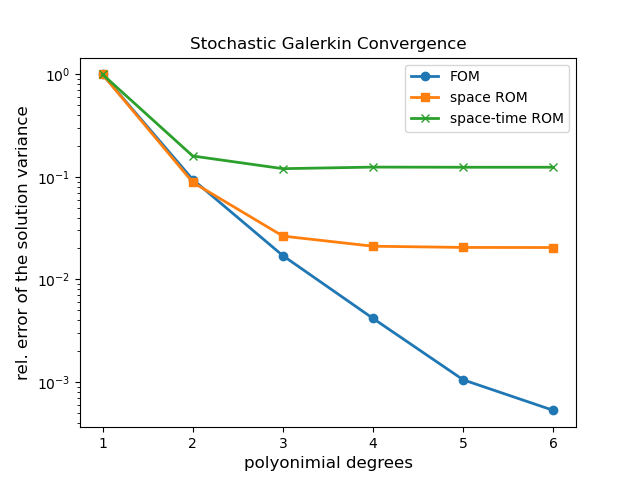}}
\caption{stochastic Galerkin convergence: relative error versus maximal polynomial degrees.}
\label{1d_sg_convergence}
\end{figure}

We discuss the convergence performances of SG FOM and ROM methods. Similarly, we plot the relative errors of solution mean and variance as quantities of interests in \Cref{1d_sg_convergence}. The relative errors in the SG setting are defined as:
\begin{equation*}
\frac{\left\|\expectation[\vState] - \expectation[\basis\reducedState]\right\|_2}{\left\|\expectation[\vState]\right\|_2}, \quad \frac{\left\|\variance[\vState] - \variance[\basis\reducedState]\right\|_2}{\left\|\variance[\vState]\right\|_2},
\end{equation*}
where the expectation and variance of the true solution obtained by 12 millions FOM MC simulations. We use integrations like \eqref{sg expectation} to calculate the expectation and variance for approximated solution function $\basis\reducedState.$ 

 \Cref{1d_sg_convergence} shows that FOM, space ROM and space-time ROM implemented by the stochastic Galerkin strategy all converge smoothly as the approximation polynomial degree increases. Similar to the MC case, space-time SG ROM has slightly larger errors than the other two methods. We again emphasize that the accuracy of the space and space-time ROM can be improved by including more basis vectors. 

\section{Numerical experiment II: 2D parametrized advection-diffusion \\problem}
In this section, we consider a rather similar advection-diffusion problem to \eqref{1d adv-diff}, but in higher dimensions and with more parameters. Specifically, we consider the following two-dimensional parametrized system with initial and boundary conditions:
\begin{equation}
\label{2d adv-diff}
\left\{
\begin{array}{l}
\displaystyle \frac{\partial \state}{\partial t} + b\, \cos(\frac{\pi}{3}) \frac{\partial \state}{\partial x} + b\, \sin(\frac{\pi}{3}) \frac{\partial \state}{\partial y} + \sigma\, \state = \nu\, \left(\frac{\partial^2 u}{\partial x^2}+\frac{\partial^2 u}{\partial y^2}\right),\\\\
\state((0,y), t, (b, \sigma, \nu)) = 0,~~\forall~y \in [0, 1], t \in [0,2.5], (b, \sigma, \nu) \in \paramDomain\\\\
\state((x,0), t, (b, \sigma, \nu)) = 0,~~\forall~x \in [0, 1], t \in [0,2.5], (b, \sigma, \nu) \in \paramDomain\\\\
\state((x, y), 0, (b, \sigma, \nu)) = 0,~~\forall~(x,y) \in [0,1] \times [0,1], (b, \sigma, \nu) \in \paramDomain
\end{array}
\right.
\end{equation}
where the state $\state: [0,1] \times [0,1] \times \paramDomain \to \RR{}$. Here, the speed $b \sim \mathcal{N}(0.5, 0.1),$ the reaction coefficient $\sigma \sim U[0.003, 0.005]$ and the diffusion coefficient $\nu \sim U[0.9, 1.1]$.

Similarly, we apply the second-order backward difference scheme for spatial discretization with $N_{s_x} = 63, N_{s_y} = 63$ nodes respectively in the $x$ and $y$ directions, hence total $N_s = N_{s_x} \times N_{s_y} = 3969$ degrees of freedom on the spatial domain. We still apply the Crank-Nicolson method with uniform time step $\Delta t = 0.005$ and implicitly solve results for all $N_t = 2.5/\Delta t = 500$ time instances. 

We set the number of training samples to be 20 in the offline phase and select a trial subspace that captures at least $99.9999\%$ of the energy in the original snapshot solution set. 

We first show the computational efficiencies of the MC space ROM and space-time ROM methods. The computation time of a ROM method is calculated as same as in \eqref{time count}.

In order to identify a spatial trial subspace for the numerical solution of \eqref{2d adv-diff}, one needs to do a singular value decomposition on a $N_s \times N_t N_{\text{train}} = 3969\times 10000$ snapshot matrix that has a large number of columns. This procedure is so expensive that it may even exceed the actual computation time of solving the reduced system. Therefore, we propose the random range finder (RRF) method \cite{HMT11} for the spatial subspace finding to reduce this overhead complexity. The methodology of RRF is introduced as follows:
\begin{enumerate}
\item
fix the number of truncated columns $\hat{K}$;
\item
right multiply a Gaussian testing matrix $\mathbf{G} \sim \mathcal{N}(0,1)^{N_t N_{\text{train}} \times \hat{K}}$ on the training solutions $\mState_{\text{train}} \in \RR{N_s \times  N_t N_{\text{train}}}$;
\item
compute the svd for the resulting matrix with a much reduced number of columns $\mState_{\text{train}} \mathbf{G} \in \RR{N_s \times \hat{K}}$ and keep the left singular vectors to be the basis.
\end{enumerate}

\subsection{Numerical results}
In \Cref{2d_mc_rom_time}, we test FOM and ROM methods on $2000$ MC samples and provide a time report for space and space-time ROMs. For the setup of RRF, we pick the number of truncated columns to be $\hat{K} = 20$ versus $K_{\text{true}} = 12$ applying a conventional svd on $\mState_{\text{train}}$. 

Some important observations from \Cref{2d_mc_rom_time} are: (1) The solving time (the third table column) of the space-time ROM method is $300$ times lower than the time of space ROM; (2) Applying RRF for spatial trial subspace reduces almost $1000$ times of the computation time compared to the conventional spatial subspace finding --- see the first table column.

\begin{table}[tbhp]
{\footnotesize
  \caption{{\bf Computation time for each step in the workflow (unit: second).}}
\label{2d_mc_rom_time}
\begin{center}
  \begin{tabular}{cccc} 
  \toprule[1pt]
 & \bf find trial subspace & \bf build ROM & \bf solve ROM \\[1mm]
\bf space ROM & $\bf \color{red} 77.908$ & $0.006$ & $63.059$\\[1.5mm]
\bf space ROM (RRF)& $\bf \color{blue} 0.822$ & $0.007$ & $76.351$\\[1.5mm]
\bf space-time ROM  & 4.016&$1.055$ & $\bf \color{blue} 0.227$\\ 
\bottomrule[1pt]
  \end{tabular}
\end{center}
}
\end{table}

\begin{figure}[tbhp]
\centering
\includegraphics[width=0.65\textwidth]{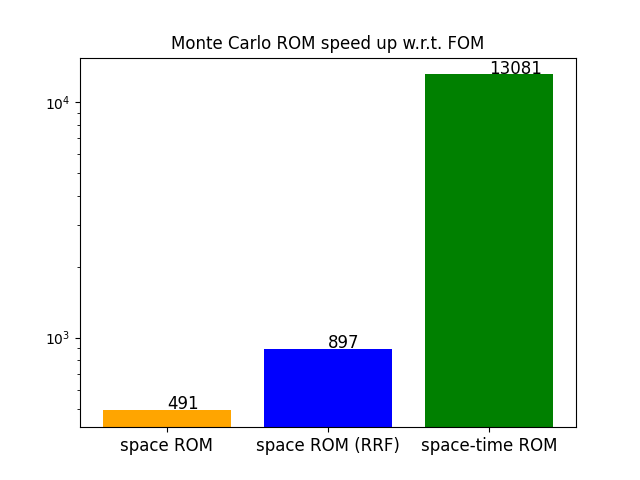}
\caption{Monte Carlo ROM speed-up.}
\label{2d_mc_speedup}
\end{figure}

From the same experiment for \Cref{2d_mc_rom_time}, we report the ROM speed-ups in \Cref{2d_mc_speedup}. The speed-up's metric is the same as in \Cref{1d_mc_speedup}. We again observe that the space-time ROM method is the fastest among the three methods. 

\begin{figure}[tbhp]
\centering
\subfigure[solution mean]{\includegraphics[width=0.48\textwidth]{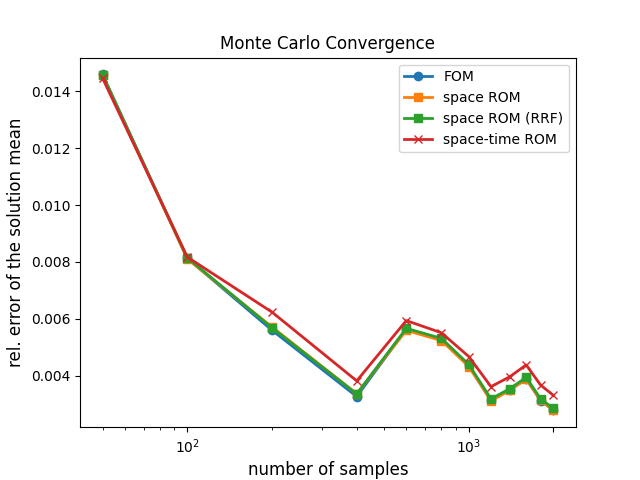}}~~
\subfigure[solution variance]{\includegraphics[width=0.48\textwidth]{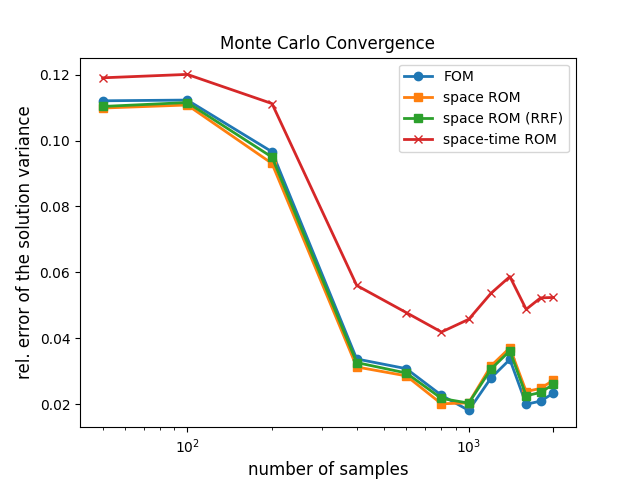}}
\caption{Monte Carlo convergence: relative error versus number of samples.}
\label{2d_mc_convergence}
\end{figure}

We next discuss the convergence property. \Cref{2d_mc_convergence} shows the convergence of the FOM and ROM methods with respect to the number of MC samples. Similar to the presentation in \Cref{1d_mc_convergence} for the 1D UQ problem, we make the y-axis to be the solution errors for relative mean and variance with respect to the number of MC samples in the x-axis. The experiment runs on a total $2000$ random instances of $(b, \sigma, \nu) \in \paramDomain$ in eq.~\eqref{2d adv-diff}. The errors are averaged results from $3$ repetitions of the test based on different random seeds.

In general, all four MC solving methods converge with the same trend. 
MC space-time ROM is the least accurate method and has the largest error means and variances. This can be considered as a trade-off of its computational efficiency and accuracy. 
The error results of FOM and both two variants of space ROM methods are very close and almost overlapping, which proves the high accuracy of MC space ROM. 

We employ stochastic Galerkin ROM methods for the 2D advection-diffusion problem \eqref{2d adv-diff}. Please note that it is commonly recognized that running a full-order model with stochastic Galerkin is too time consuming, especially in large-scale problems, thus we skip its implementation in this subsection.

\begin{figure}[tbhp]
\centering
\includegraphics[width=0.65\textwidth]{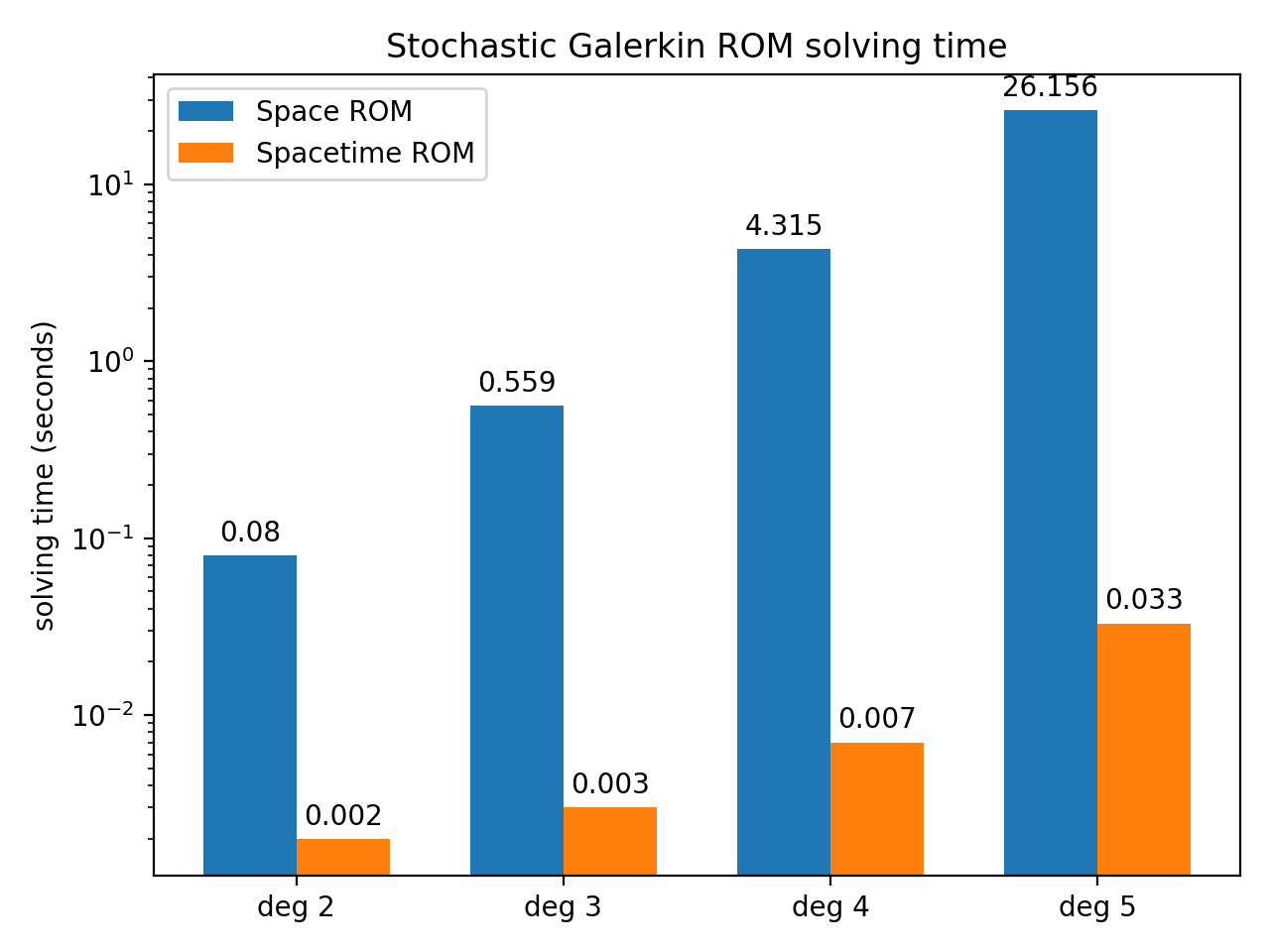}
\caption{stochastic Galerkin ROM solving time.}
\label{2d_sg_time}
\end{figure}

We are interested in the computation time of space and space-time ROMs implemented with the stochastic Galerkin strategy. In \Cref{2d_sg_time}, we report the ROM solving times (the last step in the workflow \Cref{workflow}) --- other workflow steps are already studied above in the MC case. We use the space ROM implemented by conventional subspace finding as the representative space ROM approach. 

\Cref{2d_sg_time} demonstrates the time comparisons of the SG space and space-time ROM. It is clearly shown that SG space-time ROM method is roughly $40 - 800$ times faster than SG space ROM. As the approximating polynomial degrees grows, this solving time discrepancy becomes more pronounced. 

\begin{figure}[tbhp]
\centering
\subfigure[solution mean]{\includegraphics[width=0.48\textwidth]{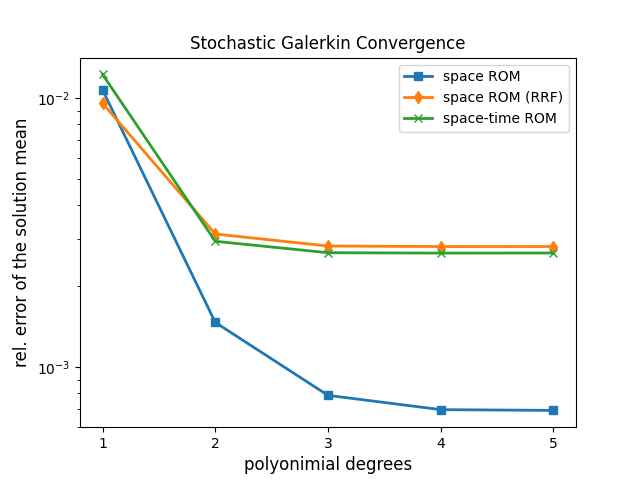}}~~
\subfigure[solution variance]{\includegraphics[width=0.48\textwidth]{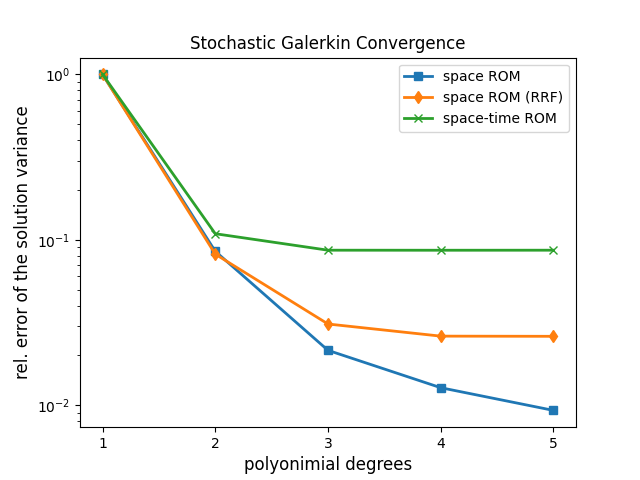}}
\caption{stochastic Galerkin convergence: relative error versus maximal polynomial degrees.}
\label{2d_sg_convergence}
\end{figure}

We finally discuss the convergence and accuracy properties of the SG ROMs. From \Cref{2d_sg_convergence}, we can observe that space ROM, space ROM with random range finder and space-time ROM methods all converge very smoothly as the approximating polynomial degree tends to increase. In terms of accuracy, the SG space ROM method with conventional trial subspace finding achieves the highest precision, while the SG space-time ROM overall has bigger error means and variances than space ROMs. 

\section{ Conclusions and future directions}
In this work, we have studied and showed the significant computational advantages of the space-time ROM method incorporated with Monte Carlo and stochastic Galerkin techniques.
By testing our proposed method on parametrized 1D and 2D advection-diffusion problems, we provided thorough numerical experiments to demonstrate both computational and convergence properties of the space-time ROM method and compared them with the FOM and space ROM methods.
The numerical performance showed that the space-time ROM method achieved remarkably high efficiency compared to the other two approaches. However, it also suffered a small loss of solution accuracy as suggested in convergence plots. 

We finally lay out the remaining future works.
\begin{enumerate}
\item
We hope to directly compare the computational cost and theoretical performance of SG ROMs to MC ROMs.
\item
Since stochastic Galerkin is an intrusive method, implementing the space-time ROM by other non-intrusive UQ propagation methods such as stochastic collocation could be a promising direction.
\item
We also want to extend the proposed method to challenging nonlinear systems where explicitly forming reduced operators is more difficult. 
\item
Developing theoretical error bounds for the proposed method is another sound future plan.
\item
It is interesting to explore the advanced sparse grid strategy to further reduce the computational complexity.
\end{enumerate}

\section{Acknowledgements}
R.\ Jin acknowledges an appointment to the NSF \\Mathematical Sciences Graduate Internship. 
This work was partially sponsored by Sandia's Advanced Simulation and
Computing (ASC) Verification and Validation (V\&V) Project/Task
\#103723/05.30.02.  This paper describes objective technical results and
analysis. Any subjective views or opinions that might be expressed in the
paper do not necessarily represent the views of the U.S.\ Department of Energy
or the United States Government.  Sandia National Laboratories is a
multimission laboratory managed and operated by National Technology and
Engineering Solutions of Sandia, LLC., a wholly owned subsidiary of Honeywell
International, Inc., for the U.S.\  Department of Energy's National Nuclear
Security Administration under contract DE-NA-0003525.

\bibliographystyle{siam}
\bibliography{references}


\end{document}